\documentclass[a4paper,12pt]{article}
\usepackage{amsthm,amsmath}
\usepackage{amssymb}
\usepackage[all]{xypic}
\usepackage{enumerate}
\usepackage{a4wide}
\newtheorem{theorem}{Theorem}[section]
\newtheorem{proposition}[theorem]{Proposition}
\newtheorem{corollary}[theorem]{Corollary}
\newtheorem{lemma}[theorem]{Lemma}
\newtheorem*{theorem*}{Theorem}
\newtheorem*{proposition*}{Proposition}
\newtheorem*{corollary*}{Corollary}
\newtheorem*{lemma*}{Lemma}\theoremstyle{definition}
\newtheorem{definition}[theorem]{Definition}

\newtheorem{example}[theorem]{Example}
\newtheorem{remark}[theorem]{Remark}

\newtheorem*{definition*}{Definition}

\newcommand{\cat}[1]{\mathcal{#1}}
\newcommand{\coring}[1]{\mathfrak{#1}}
\newcommand{\tensor}[1]{\otimes_{#1}}

\newcommand{\rcomod}[1]{\mathcal{M}^{#1}}
\newcommand{\lcomod}[1]{{}^{#1}\mathcal{M}}
\newcommand{\rmod}[1]{\mathcal{M}_{#1}}
\newcommand{\lmod}[1]{{}_{#1}\mathcal{M}}
\newcommand{\bimod}[2]{{}_{#1}\mathcal{M}_{#2}}
\newcommand{\bcomod}[2]{{}^{#1}\mathcal{M}^{#2}}

\renewcommand{\hom}[3]{\mathrm{Hom}_{#1}(#2,#3)}
\newcommand{\End}[1]{\mathrm{End}(#1)}
\newcommand{\rend}[2]{\mathrm{End}({#2}_{#1})}
\newcommand{\lend}[2]{\mathrm{End}({}_{#1}{#2})}

\newcommand{\esc}[2]{\langle #1,#2 \rangle}
\newcommand{\Rat}{\mathrm{Rat}}

\newcommand{\naturales}{\mathbb{N}}

\newcommand{\rring}{\coring{C}^{\ast}}
\title{Semisimple corings}

\author{L. El Kaoutit \and J. G\'omez-Torrecillas \and F. J. Lobillo}

\date{}

\begin{document}
\maketitle

\begin{abstract}
This paper states the basic essentials for a theory of semisimple
corings.
\end{abstract}

\section*{Introduction}
M. Sweedler \cite{Sweedler:1975} introduced the notion of coring
as a generalization of the concept of coalgebra in order to study
the set of intermediate division rings for an extension of
division rings. It turns out that this formalism embodies several
kinds of relative module categories. Thus, graded modules,
Doi-Hopf modules and, more generally, entwined modules are
instances of comodules over suitable corings (see
\cite{Brzezinski:2000unp} and its references). From this point of
view, an interesting question is to characterize those corings
encoding the simplest type of category of comodules. Since, under
certain hypotheses, the category of comodules is abelian, the
simple objects play a relevant role in its structure. In the most
favorable case all comodules are direct sums of simple comodules
(that is, the category of comodules is semisimple). In the
classical theory of modules over rings, the study of semisimple
rings precedes the development of the entire theory. This paper
states the basic essentials for a theory of semisimple corings.

Throughout this paper the word ring will refer to an associative
algebra over
 a commutative ring $K$, and
the term subring is then understood as subalgebra. The category
of all left modules over a ring $R$ will be denoted by
$\lmod{R}$, being $\rmod{R}$ the notation for the category of all
right $R$--modules. An agile introduction to abelian categories is
contained in \cite{Stenstrom:1975}. The notation $X \in \cat{A}$
for a category $\cat{A}$ means that $X$ is an object of
$\cat{A}$, and the identity morphism attached to any object $X$
will be represented by the object itself.

\section{Corings and comodules}

We discuss under which conditions the category of (right) comodules over a coring is abelian.
We first recall from \cite{Sweedler:1975} the notion of a coring.

Let $A$ denote a ring. An $A$--\emph{coring} is a three-tuple
($\coring{C}, \Delta,\epsilon$) consisting of an $A$-bimodule
$\coring{C}$ and two $A$--bimodule maps
\[
\Delta : \coring{C} \longrightarrow \coring{C} \tensor{A}
\coring{C}, \hspace{2em} \epsilon : \coring{C} \longrightarrow A
\]
such that the diagrams
\begin{displaymath}
\xymatrix@C=60pt{ \coring{C} \ar^-{\Delta}[r] \ar_-{\Delta}[d] &
\coring{C} \tensor{A} \coring{C} \ar^-{\coring{C} \tensor{A}
\Delta}[d] \\ \coring{C} \tensor{A} \coring{C} \ar^-{ \Delta
\tensor{A} \coring{C}}[r] & \coring{C} \tensor{A} \coring{C}
\tensor{A} \coring{C} }
\end{displaymath}
and
\begin{displaymath}
\xymatrix{\coring{C} \ar^-{\Delta}[r] \ar_-{\cong}[dr] &
\coring{C} \tensor{A} \coring{C} \ar^-{\coring{C} \tensor{A}
\epsilon}[d] \\ & \coring{C} \tensor{A} A } \hspace{3em}
\xymatrix{\coring{C} \ar^-{\Delta}[r] \ar_-{\cong}[dr] &
\coring{C} \tensor{A} \coring{C} \ar^-{\epsilon \tensor{A}
\coring{C}}[d] \\ & A \tensor{A} \coring{C} }
\end{displaymath}
commute.

A \emph{left} $\coring{C}$--\emph{comodule} is a pair ($M,
\lambda_M$) consisting of a left $A$--module $M$ and an $A$-linear
map $\lambda_M : M \rightarrow \coring{C} \tensor{A} M$ such that
the diagrams
\begin{displaymath} \xymatrix@C=60pt{ M
\ar^-{\lambda_M}[r] \ar_-{\lambda_M}[d] & \coring{C} \tensor{A} M
\ar^-{\Delta \tensor{A} M}[d] \\ \coring{C} \tensor{A} M \ar^-{
\coring{C} \tensor{A} \lambda_M}[r] & \coring{C} \tensor{A}
\coring{C} \tensor{A} M } \hspace{3em} \xymatrix{ M
\ar^-{\lambda_M}[r] \ar_-{\cong}[dr] & \coring{C} \tensor{A} M
\ar^-{\epsilon \tensor{A} M}[d] \\ & A \tensor{A} M}
\end{displaymath}
commute. Right $\coring{C}$--comodules are
similarly defined; we use the notation $\rho_M$ for their
structure maps. A \emph{morphism} of left $\coring{C}$--comodules
($M,\lambda_M$) and ($N,\lambda_N$) is an $A$-linear map $f: M
\rightarrow N$ such that the following diagram is commutative
\begin{displaymath}
\xymatrix@C=60pt{ M \ar^-{f}[r]
\ar_-{\lambda_M}[d] & N \ar^-{\lambda_N}[d] \\ \coring{C}
\tensor{A} M \ar^-{\coring{C} \tensor{A} f}[r] & \coring{C}
\tensor{A} N. }
\end{displaymath}
The $K$--module of all left $\coring{C}$--comodule morphisms from
$M$ to $N$ is denoted by $\hom{\coring{C}}{M}{N}$. The category
of all left $\coring{C}$--comodules will be denoted by
$\lcomod{\coring{C}}$. Analogously, we can consider the category
of all right $\coring{C}$--comodules $\rcomod{\coring{C}}$. Every
valid statement about left comodules entails a correct assertion
for right comodules, which will be implicitly understood.

Coproducts and cokernels in $\lcomod{\coring{C}}$ do exist, and
they can be already computed in $\lmod{A}$. Therefore,
$\lcomod{\coring{C}}$ has arbitrary inductive limits. If
$\coring{C}_A$ is a flat module, then $\lcomod{\coring{C}}$ is
easily proved to be an abelian category. The converse is not
true, as the following example shows.

\begin{example}\label{elejemplo}
Let $A = \begin{pmatrix}
  R & B \\
  0 & S
\end{pmatrix}$ be a generalized triangular matrix ring with $B$ an
$(R,S)$--bimodule over the rings $R$ and $S$. It is well-known
that the right $A$--modules are given by three-tuples $M =
(M',M'',\mu)$ consisting of a right $R$--module $M'$, a right
$S$--module $M''$ and an $S$--module map $\mu : M' \tensor{R} B
\rightarrow M''$. A homomorphism of right $A$--modules is then
given by a pair $f = (f',f'') : (M',M'',\mu) \rightarrow
(N',N'',\nu)$ consisting of an $R$--module map $f' : M'
\rightarrow N'$ and a right $S$--module map $f'' : M'' \rightarrow
N''$ such that $f''\mu = \nu (f' \tensor{R} B)$. Consider the
ideal $I =
\begin{pmatrix}
  R & B \\
  0 & 0
\end{pmatrix}$ of $A$ which, as a right $A$--module, corresponds
to $(R,B,\mu)$, where $\mu : R \tensor{R} B \rightarrow B$ is the
canonical isomorphism. Now, $M \tensor{A} I = (M', M'\tensor{R}
B, id)$ and the multiplication map $M \tensor{A} I \rightarrow M$
is given by $(id,\mu) : (M', M' \tensor{R} B, id) \rightarrow
(M',M'',\mu)$. The $A$--bimodule $I$ is an $A$--coring with
comultiplication given by the isomorphism $I \cong I \tensor{A}
I$, and counit given by the inclusion $I \subseteq A$. It can be
easily shown that a right $A$--linear coaction $\rho_M =
(\rho',\rho'') : M \rightarrow M \tensor{A} I$ is an
$I$--comodule structure if and only if $\rho' = id_{M'}$ and
$\mu$ is an isomorphism with inverse $\rho''$. Therefore, the
category $\rcomod{I}$ of all right $I$--comodules can be
identified with the category of all right $A$--modules
$(M',M'',\mu)$ such that $\mu$ is an isomorphism. Now the functor
$F : \rcomod{I} \rightarrow \rmod{R}$ given by $F(M',M'',\mu) =
M'$ is easily shown to be an equivalence of categories. In
particular, $\rcomod{I}$ is a Grothendieck category and ${}_AI$
is not flat unless ${}_RB$ is.
\end{example}

The following result clarifies the situation created by our
example. Recall from \cite[Proposition 3.1]{Guzman:1989} that the
functor $\coring{C}\tensor{A}- \,: \lmod{A} \rightarrow
\lcomod{\coring{C}}$ is right adjoint to the forgetful functor $U
\,: \lcomod{\coring{C}} \rightarrow \lmod{A}$.

\begin{proposition}\label{plano}
Let $\coring{C}$ be an $A$--coring. Then the following are
equivalent
  \begin{enumerate}[(i)]
  \item $\coring{C}_{A}$ is flat.
  \item $\lcomod{\coring{C}}$ is an abelian category and
the functor $U$ is left exact.
    \item $\lcomod{\coring{C}}$ is a Grothendieck category and
the functor $U$ is left exact.
  \end{enumerate}
  \begin{proof}
  $(i) \Rightarrow (iii)$ The exactness of the functor
  $\coring{C} \tensor{A} - : \lmod{A} \rightarrow \lmod{A}$ entails that kernels in
  $\lcomod{\coring{C}}$ can be already computed in $\lmod{A}$. This gives that
  $\lcomod{\coring{C}}$ is
  a complete and co-complete abelian category with exact direct limits. We need to find a generator
  for $\lcomod{\coring{C}}$. For this, we proceed as in the
  proof of \cite[13.13]{Wisbauer:1998}. Let $M
\in \lcomod{\coring{C}}$, with coaction $\lambda_M$, and $A^{(I)}
\rightarrow M \rightarrow 0$ the free presentation of $M$ in
$\lmod{A}$, we have
\begin{displaymath}
\xymatrix{ \coring{C}^{(I)} \cong \coring{C} \tensor{A} A^{(I)}
\ar^-{g}[r] & \coring{C} \tensor{A} M \ar[r] & 0 \\  g ^{-1}(M)
\ar[r] \ar@{^{(}->}[u] & M \ar@{^{(}->}_-{\lambda_M}[u] \ar[r] & 0
}
\end{displaymath}
now it is clear that
\[
 \mathfrak{g}^{l} = \oplus \{ L |\,L \text{ is left subcomodule of }
\coring{C}^k, \,k \in \naturales\}
\]
is a generator of $\lcomod{\coring{C}}$. Obviously, the forgetful functor
$U : \lcomod{\coring{C}} \rightarrow \lmod{A}$ is exact in this case. \\
$(iii) \Rightarrow (ii)$ This is obvious. \\
  $(ii) \Rightarrow (i)$ By \cite[Corollary 3.2.3]{Popescu:1973},  $\coring{C}
\tensor{A}- \,: \lmod{A} \rightarrow \lcomod{\coring{C}}$ is left
exact and, thus, $U \circ (\coring{C} \tensor{A}-) \,: \lmod{A}
\rightarrow \lmod{A}$ is a left exact too. Therefore,
$\coring{C}_{A}$ is a flat module.
  \end{proof}
\end{proposition}

As a consequence of the proof of Proposition \ref{plano} we obtain:

\begin{corollary}\label{Csubgenera}
If $\coring{C}_A$ is flat then every left $\coring{C}$--comodule is
isomorphic to a subcomodule
of a $\coring{C}$--generated left comodule.
\end{corollary}

\begin{example}\label{ejemplo1}(\cite[Example 1.2]{Sweedler:1975}). Let $B
\rightarrow A$ a ring extension. Then $\coring{C} = A \tensor{B}
A$ is an $A$--coring with the coproduct $\coring{C} \rightarrow
\coring{C} \tensor{A} \coring{C}$ given by $a \tensor{B} a'
\mapsto a \tensor{B} 1 \tensor{A} 1 \tensor{B} a'$ and the counit
$A \tensor{B} A \rightarrow A$ is the multiplication map.
\end{example}

\begin{example}(\cite[Proposition 2.2]{Brzezinski:2000unp}) Let
$(C,\Delta,\epsilon)$ be a $K$--coalgebra and assume that the
canonical left $A$--module $\coring{C} = A \tensor{} C$ has a
right $A$--module structure which makes $\coring{C}$ an
$A$--bimodule. Define the $K$--linear map $\psi: C \tensor{} A
\rightarrow A \tensor{} C$ as $\psi(c \tensor{} a) = (1 \tensor{}
c)a$. Consider the left $A$--module maps $\Delta_{\coring{C}}:
\coring{C} \rightarrow \coring{C} \tensor{A} \coring{C} \cong A
\tensor{} C \tensor{} C$, $\Delta_{\coring{C}} = A \tensor{}
\Delta$, and $\epsilon_{\coring{C}} = A \tensor{} \epsilon$. Then
$(\coring{C},\Delta_{\coring{C}},\epsilon_{\coring{C}})$ is an
$A$--coring if and only if $(A,C)_{\psi}$ is an entwining
structure (see \cite{Brzezinski/Majid:1998}). Moreover, the
category of comodules $\rcomod{\coring{C}}$ is isomorphic to the
category of entwined modules $\mathcal{M}_A^C(\psi)$. Examples of
categories of entwined modules are Doi-Koppinen modules,
introduced in \cite{Doi:1992} and \cite{Koppinen:1994} (cf.
\cite[Example 3.1(3)]{Brzezinski:1999}).
\end{example}

\section{Rational modules and comodules}

We state a formal framework (the notion of rational pairing) which
reduces the study of some categories of comodules to the
investigation of certain subcategories of categories of modules.
The development is adapted from the given in \cite{Gomez:1998} and
\cite{Abuhlail/Gomez/Lobillo:2001} for coalgebras over commutative
rings. A somewhat different approach is \cite{Wisbauer:2001unp}.

Let $P,Q$ be
$A$--bimodules. Any balanced bilinear form \begin{displaymath}
\esc{-}{-} : P \times Q \longrightarrow A \end{displaymath}
provides natural transformations $\alpha : Q \tensor{A} -
\rightarrow \hom{A}{{}_AP}{-}$ and $\beta : - \tensor{A} P
\rightarrow \hom{A}{Q_A}{-}$ given by \begin{gather*}
\begin{split} \alpha_M : Q \tensor{A} M &\longrightarrow
\hom{A}{_{A}P}{{}_{A}M} \\ q \tensor{A} m &\longmapsto \left[ p
\mapsto \esc{p}{q}m \right] \end{split} \\ \begin{split} \beta_N
: N \tensor{A} P &\longrightarrow \hom{A}{Q_{A}}{N_{A}} \\ n
\tensor{A} p &\longmapsto \left[ q \mapsto n\esc{p}{q} \right].
\end{split} \end{gather*} Moreover, if $M$ is an $A$--bimodule
then $\alpha_M$ and $\beta_M$ are bimodule homomorphisms. The
canonical isomorphisms provide two bimodule maps \begin{gather}
\begin{split} \alpha_A : Q &\longrightarrow
\hom{A}{_{A}P}{{}_{A}A} = {}^*P \\ q &\longmapsto \left[ p
\mapsto \esc{p}{q} \right] \end{split} \label{eq:alpha} \\
\begin{split} \beta_A : P &\longrightarrow \hom{A}{Q_A}{A_A}= Q^*
\\ p &\longmapsto \left[ q \mapsto \esc{p}{q} \right] \end{split}
\label{eq:beta} \end{gather} which are bimodule homomorphisms. So
we can recover the balanced bilinear form if one of the natural
transformations is given.

\begin{definition}
The data $T = (P,Q,\esc{-}{-})$ are called a \emph{left rational
system} if $\alpha_M$ is injective for each left $A$--module $M$,
and a \emph{right rational system} if $\beta_N$ is injective for
every right $A$--module $N$.
\end{definition}

\begin{remark}\label{plano:rational}
(\cite[Remark 2.4]{Abuhlail/Gomez/Lobillo:2001}). Let
$(P,Q,\esc{-}{-})$ be a left rational system. Let $M \in \lmod{A}$
and $N$ be submodule of $M$ with the canonical injection $i_N$.
Consider the following commutative diagram
\begin{displaymath}
\xymatrix{ Q \tensor{A} N \ar@{^{(}->}^-{\alpha_{N}}[r] \ar_-{Q
\tensor{A} i_N }[d] & \hom{A}{P}{N} \ar@{^{(}->}^-{i}[d]\\ Q
\tensor{A} M \ar@{^{(}->}^-{\alpha_M}[r] & \hom{A}{P}{M}. }
\end{displaymath}
Hence $Q \tensor{A} i_N$ is injective. Since $M$ was an arbitrary
left $A$--module, we conclude that $Q_A$ should be a flat
$A$--module. Analogously, if $(P,Q,\esc{-}{-})$ is a right
rational system, then we get that ${}_{A}P$ is flat.
\end{remark}

Let $\esc{-}{-}:P \times Q \rightarrow A$ and $[-,-]:P' \times Q'
\rightarrow A$ be two balanced bilinear forms with natural
transformations $\alpha,\beta$ and $\alpha',\beta'$ respectively.
We can define a new balanced bilinear form \begin{displaymath}
\begin{split} \{-,-\}:(P \tensor{A} P') \times (Q' \tensor{A} Q)
&\longrightarrow A \\ (p \tensor{A} p',q' \tensor{A} q)
&\longmapsto \{p \tensor{A} p' , q' \tensor{A} q\} =
\esc{p}{[p',q']q} = \esc{p[p',q']}{q}. \end{split}
\end{displaymath} The natural transformations associated to
$\{-,-\}$ are given by the compositions
\begin{displaymath}
\xymatrix{ Q' \tensor{A} Q \tensor{A} M \ar_-{\alpha'_{Q
\tensor{A} M}}[d] \ar@{-->}[r] & \hom{A}{P \tensor{A} P'}{M} \\
\hom{A}{_{A}P'}{Q \tensor{A} M} \ar_-{(\alpha_M)_*}[r] &
\hom{A}{_{A}P'}{\hom{A}{_{A}P}{M}}, \ar_-{i_M}[u]}
\end{displaymath}
\begin{displaymath} \xymatrix{ N \tensor{A} P
\tensor{A} P' \ar_-{\beta'_{N \tensor{A} P}}[d] \ar@{-->}[r] &
\hom{A}{Q' \tensor{A} Q}{N} \\ \hom{A}{Q'_{A}}{N \tensor{A} P}
\ar_-{(\beta_N)_*}[r] & \hom{A}{Q'_{A}}{\hom{A}{Q_{A}}{M}}.
\ar_-{i_N}[u]}
\end{displaymath}

The following proposition, which is now clear, replaces \cite[Proposition 2.2]{Gomez:1998}
in order to show that the canonical comodule structure over a rational module is
pseudocoassociative.

\begin{proposition}\label{bd:besc} Let $(P,Q,\esc{-}{-})$ and
$(P',Q',[-,-])$ be two left (resp. right) rational systems. Then
$(P \tensor{A} P',Q' \tensor{A} Q, \{-,-\})$ is also a left
(resp. right) rational system. \end{proposition}

Let $(\coring{C},\Delta,\epsilon)$ be an $A$--coring. Recall that
$\coring{C}^* = \hom{A}{\coring{C}_A}{A}$ (resp. ${}^*\coring{C}
= \hom{A}{{}_A\coring{C}}{A}$) is a ring extension of $A^{op}$
with multiplication $gf = f \circ (g \tensor{A} \coring{C}) \circ
\Delta$ (resp. $gf = g \circ (\coring{C} \tensor{A} f) \circ
\Delta$). Both units are $\epsilon$. See \cite[Proposition
3.2]{Sweedler:1975} for details.

\begin{definition} A \emph{left rational pairing} is a left
rational system $(B,\coring{C},\esc{-}{-})$ such that $B$ is a
ring extension of $A$, $\coring{C}$ is an $A$--coring, and $\beta
: B \rightarrow \coring{C}^*$ is a ring antimorphism. If
$\Delta(c) = \sum_i c_i \tensor{A} d_i$ then
\begin{equation}\label{eq:betaring} \esc{ab}{c} = \esc{a}{\sum_i
\esc{b}{c_i}d_i},\, \forall a,\, b \in B \quad \text{and} \quad
\epsilon = \esc{1}{-}. \end{equation} Analogously, a \emph{right
rational pairing} is a right rational system
$(\coring{C},B',[-,-])$ such that $B'$ is a ring extension of $A$,
$\coring{C}$ is an $A$--coring and $\alpha : B' \rightarrow
{}^*\coring{C}$ is a ring antimorphism. If $\Delta(c) = \sum_i c_i
\tensor{A} d_i$ then

\begin{displaymath}
[c,ab] = [\sum_i c_i[d_i,a],b],\, \forall a,\, b \in B \quad
\text{and} \quad \epsilon = [1,-].
\end{displaymath}
\end{definition}

\begin{example}\label{ejemplo:rat} Let $\coring{C}$ be an
$A$--coring such that $\coring{C}$ is projective as a right
$A$--module. By using any dual basis associated with the
projectivity of $\coring{C}_A$, we prove that the canonical
balanced bilinear form $\esc{-}{-}:{\coring{C}^*}^{op} \times
\coring{C} \rightarrow A$ gives a left rational pairing $T =
({\coring{C}^*}^{op},\coring{C},\esc{-}{-})$. Analogously, if
$\coring{C}$ is an $A$--coring such that ${}_{A}\coring{C}$ is a
projective module, then
$T'=(\coring{C},{}^{\ast}\coring{C}^{op},[-,-])$ is a right
rational pairing. \end{example}

Let $T = (B,\coring{C},\esc{-}{-})$ be a left rational pairing. An
element $m$ in a left $B$--module $M$ is called \emph{rational} if
there exists a set of \emph{left rational parameters}
$\{(c_i,m_i)\} \subseteq \coring{C} \times M$ such that $b m =
\sum_i \esc{b}{c_i} m_i$ for all $b \in B$. The set of rational
elements in $M$ is denoted by $\Rat^T(M)$. The proofs detailed in
\cite[Section 2]{Gomez:1998} can be adapted in a straightforward
way in order to get that $\Rat^T(M)$ is a $B$--submodule of $M$
and that the assignment $M \mapsto \Rat^T(M)$ defines a functor
\begin{equation*} \Rat^T : \lmod{B} \rightarrow \lmod{B},
\end{equation*} which is in fact a left exact preradical.
Therefore, the full subcategory $\Rat^T(\lmod{B})$ of $\lmod{B}$
whose objects are those $B$--modules $M$ such that $\Rat^T(M) =
M$ is a closed reflective subcategory \cite[p. 395]{Gabriel:1962}
and, in particular, it is a Grothendieck category. The modules in
the subcategory $\Rat^T(\lmod{B})$ will be called \emph{rational
left $B$--modules} (with respect to $T$). Now it turns out that
every rational left $B$--module is a left $\coring{C}$--comodule
with structure map $\lambda_M : M \rightarrow \coring{C}
\tensor{A} M$ given by $\lambda_M (m) = \sum c_i \tensor{A} m_i$,
where $\{ (c_i,m_i) \}$ is any set of rational parameters for $m
\in M$ (\cite[Proposition 3.5]{Gomez:1998} for a proof which can
be adapted to the present setting). This leads to a functor
\begin{equation*} {}^{\coring{C}}(-) : \Rat^T(\lmod{B})
\longrightarrow \lcomod{\coring{C}} \end{equation*} which can be
shown to be an isomorphism of categories with the guide of
\cite[Section 3]{Gomez:1998}. It can be also deduced that ${}_B
\coring{C}$ becomes a subgenerator for $\Rat^T(\coring{C})$.
Therefore, we can state

\begin{theorem}\label{thm:iso}
Let $T = (B, \coring{C}, \esc{-}{-})$ be a left rational pairing.
The functor ${}^{\coring{C}}(-) : \Rat^T(\lmod{B}) \rightarrow
\lcomod{\coring{C}}$ is an isomorphism of categories.
Moreover, every left $\coring{C}$--comodule
is isomorphic to a $B$--submodule of a
${}_B\coring{C}$-generated $B$--module.
\end{theorem}

This theorem, when applied to the rational pairing
$T = ({\coring{C}^*}^{op},\coring{C},\esc{-}{-})$
given in Example \ref{ejemplo:rat} leads to

\begin{corollary}\label{ratcom} Let $\coring{C}$ an $A$--coring.
If $\coring{C}_A$ is projective, then the functor
$(-)^{\coring{C}} : \Rat^{T}(\rmod{\coring{C}^*}) \rightarrow
\lcomod{\coring{C}}$ is an isomorphism of categories. Moreover,
every left $\coring{C}$--comodule is isomorphic to a
$\coring{C}^*$--submodule of a
$\coring{C}_{\coring{C}^*}$--generated $\coring{C}^*$--module.
\end{corollary}

Theorem \ref{thm:iso} has a right analogue. If $T' =
(\coring{C},B',[-,-])$ is a right rational pairing, then we can
define functors $(-)^{\coring{C}} : \Rat^{T'}(\rmod{B'})
\rightarrow \rcomod{\coring{C}}$ and $(-)_{B'} :
\rcomod{\coring{C}} \rightarrow \Rat^{T'}(\rmod{B'})$. These
functors lead to the following theorem

\begin{theorem*}\textbf{\upshape \ref{thm:iso}'.}
Let $T' = (\coring{C}, B',[-,-])$ be a right rational pairing. The
functor $(-)^{\coring{C}} : \Rat^{T'}(\rmod{B'}) \rightarrow
\rcomod{\coring{C}}$ is an isomorphism of categories. Moreover,
every right $\coring{C}$--comodule is isomorphic to a
$B'$--submodule of a $\coring{C}_{B'}$-generated $B'$--module.

\end{theorem*}

Finally, we state a useful consequence of the former development.

\begin{proposition}\label{f.g}
Let $T = (B, \coring{C}, \esc{-}{-})$  be a left (resp. $T'=
(\coring{C},B',[-,-])$ right) rational pairing. Let $M \in
\lcomod{\coring{C}}$. Then $M$ is a finitely generated left (resp.
right) $\coring{C}$--comodule if and only if $M$ is finitely
generated left  (resp. right) $A$--module.
\end{proposition}

Recall that a $\coring{C}$--bicomodule is an $A$--bimodule $M$
endowed with a right $A$--linear left $\coring{C}$--comodule
structure $\lambda_M : M \rightarrow \coring{C} \tensor{A} M$ and
a left $A$--linear right $\coring{C}$--comodule structure $\rho_M
: M \rightarrow M \tensor{A} \coring{C}$ such that
\begin{equation}\label{bicomcon}
(\lambda_M \tensor{A} \coring{C})\rho_M = (\coring{C} \tensor{A}
\rho_M) \lambda_M.
\end{equation}
The $\coring{C}$--bicomodules are the objects of a
category $\bcomod{\coring{C}}{\coring{C}}$
whose morphisms are those $A$--bimodule homomorphisms
which are left and right
$\coring{C}$--colinear.

Let $T = (B,\coring{C},\esc{-}{-})$
(resp. $T' = (\coring{C},B',[-,-])$) be a left
(resp. right) rational pairing.

\begin{lemma}\label{subbisubcom} Let $M$ be an $A$--bimodule with
a left $\coring{C}$--comodule structure $\lambda_M : M \rightarrow
\coring{C} \tensor{A} M$ and a right $\coring{C}$--structure map
$\rho_M : M \rightarrow M \tensor{A} \coring{C}$. Then $M$ is a
$\coring{C}$--bicomodule if and only if $M$ is a
$(B,B')$--bimodule. \end{lemma} \begin{proof} By Theorems
\ref{thm:iso} and \ref{thm:iso}', $M$ is a rational left
$B$--module and a rational right $B'$--module. We first prove that
$\lambda_M$ is right $A$--linear if and only if $M$ is a
$(B,A)$--bimodule as follows: for each $m \in M$, write $\lambda_M
(m) = \sum c_i \tensor{A} m_i$ for a set of left rational
parameters $\{ (c_i, m_i) \}$. Thus, $\lambda_M$ is right
$A$--linear if and only if $\{(c_i,m_ia) \}$ is a set of rational
parameters for $ma$ for our generic $m$ and every $a \in A$. But
this last condition is easily proved to be equivalent to require
that $M$ is a $(B,A)$--bimodule. Of course, $\rho_M$ is left
$A$--linear if and only if $M$ is an $(A,B')$--bimodule. Thus we
see that, in order to prove the Lemma, we can assume that $M$ is a
$(B,A)$--bimodule and an $(A,B')$--bimodule. Under this condition,
$M$ is a $\coring{C}$--bicomodule if and only if
\[
\sum c_i \tensor{A} m_{ij} \tensor{A} d_{ij} = \sum c_{ji}
\tensor{A} m_{ji} \tensor{A} d_{j},
\]
where $\{(c_i,m_i)\}$ is a set of left rational parameters for
$m$, $\{(m_{ij},d_{ij})\}$ is a set of right rational parameters
for each $m_i$, $\{(m_j,d_j)\}$ is a set of right rational
parameters for $m$, and $\{(c_{ji},m_{ji})\}$ is a set of left
rational parameters for each $m_j$. An easy computation gives
\begin{displaymath} \begin{array}{lll} (b.m).b' & =\sum
\esc{b}{c_i}m_{ij}[d_{ij},b'], &  \text{ and} \\ b.(m.b')& = \sum
\esc{b}{c_{ji}}m_{ji}[d_j,b'], & \text{ for any } (b,b') \in B
\times B'. \end{array} \end{displaymath} Hence $M$ is
$(B,B')$--bimodule if and only if \begin{displaymath} \sum
\esc{b}{c_i}m_{ij}[d_{ij},b']=\sum \esc{b}{c_{ji}}m_{ji}[d_j,b'],
\text{ for any } (b,b') \in B \times B'.
\end{displaymath}
By Remark \ref{plano:rational}, the following map is injective
\begin{displaymath}
\xymatrix@1@=45pt{\coring{C} \tensor{A} M \tensor{A} \coring{C}
\ar@{^{(}->}^{\coring{C} \tensor{A} \beta_M }[r] & \coring{C}
\tensor{A} \hom{A}{B'_{A}}{M}
\ar@{^{(}->}^{\alpha_{\hom{A}{B'_{A}}{M}}}[r] &
\hom{A}{{}_{A}B}{\hom{A}{B'_A}{M}} }.
\end{displaymath}
Hence $M$ is a $\coring{C}$--bicomodule if and only it is a
$(B,B')$--bimodule.
\end{proof}

\begin{proposition}\label{bi-rational} Let $T=
(B,\coring{C},\esc{-}{-})$ be and $T'= (\coring{C},B',[-,-])$ be
rational pairings, and let $\Rat^{T,T'}(\bimod{B}{B'})$ the full
subcategory of the category $\bimod{B}{B'}$ whose objects are the
$(B,B')$--bimodules which are rational as $B$--modules and as
$B'$--modules. Then there is an isomorphism of categories
$\Rat^{T,T'}(\bimod{B}{B'}) \cong
\bcomod{\coring{C}}{\coring{C}}$. \end{proposition} \begin{proof}
If $M$ is a $\coring{C}$--bicomodule, then, by Lemma
\ref{subbisubcom}, $M$ is a $(B,B')$--bimodule and, by Theorems
\ref{thm:iso} and \ref{thm:iso}~', $M$ is rational as a left
$B$--module and as a right $B$--module. Conversely, every
$(B,B')$--bimodule $M$ such that the modules ${}_BM$ and $M_{B'}$
are rational is, by Lemma \ref{subbisubcom} and theorems
\ref{thm:iso} and \ref{thm:iso}', a $\coring{C}$--bicomodule.
\end{proof}

\begin{corollary}\label{subcosubbi}
Let $\coring{I}$ be an $A$--sub-bimodule of $\coring{C}$.
\begin{enumerate}
\item $\coring{I}$ is a sub-bicomodule of $\coring{C}$ if and only if $\coring{I}$ is a
$(B,B')$--sub-bimodule of $\coring{C}$.
\item If $\coring{I}$ is
pure both as a left and a right $A$--submodule of $\coring{C}$,
then $\coring{I}$ is a subcoring of $\coring{C}$ if and only if it
is a $(B,B')$--sub-bimodule.
\end{enumerate}
\end{corollary}

For a left $\coring{C}$--comodule $M$ define $C(M)$ as the sum of
the images of all comodule homomorphisms from $M$ to $\coring{C}$.
In presence of left and right rational pairings $T =
(B,\coring{C},\esc{-}{-})$ and $T' = (\coring{C},B',[-,-])$, it is
easy to prove that $C(M)$ is a sub-bicomodule of $\coring{C}$: by
definition, it is a left $B$--submodule of the $(B,B')$--bimodule
$\coring{C}$. Now, if $b' \in B'$ and $c = f(m)$ for some
homomorphism of left $\coring{C}$--comodules, then $cb' = (r_{b'}
\circ  f)(m)$, where $r_{b'} : \coring{C} \rightarrow \coring{C}$
is the homomorphism of left comodules given by right
multiplication by $b'$. Therefore $C(M)$ is a
$(B,B')$--sub-bimodule of $\coring{C}$ and, by Corollary
\ref{subcosubbi}, it is a sub-bicomodule of $\coring{C}$, which
will be called \emph{bicomodule of coefficients} of $M$.

\begin{proposition}\label{C(M)} Let $\lambda_M : M \rightarrow
\coring{C} \tensor{A} M$ be a left comodule and assume there are
rational pairings $T = (B,\coring{C},\esc{-}{-})$ and $T' =
(\coring{C},B',[-,-])$ on the left and on the right,
respectively.
\begin{enumerate}
\item If $\tau :\coring{I}\hookrightarrow \coring{C}$ is a
monomorphism of $\coring{C}$--bicomodule, such that $\lambda_M(M)
\subseteq (\tau \tensor{A} M)(\coring{I} \tensor{A} M)$, then
$C(M) \subseteq \coring{I}$.
\item If $N$ is a subcomodule of $M$, then $C(N) \subseteq C(M)$
and $C(M/N) \subseteq C(M)$. \item If $N \cong M$ is an
isomorphism of comodules, then $C(N) = C(M)$.
\end{enumerate}
\end{proposition}
\begin{proof} Let $c = f(m) \in
C(M)$, where $f : M \rightarrow \coring{C}$ is a homomorphism of
left comodules, and write $\lambda_M (m) = \sum c_i \tensor{}
m_i$, for some $c_i \in \coring{J}$ and $m_i \in M$. Since $f$ is
a comodule map, we have
\[
\Delta(c) = \Delta(f(m)) = (\coring{C} \tensor{A} f)(\lambda_{M})(m) =
\sum c_i \tensor{} f(m_i),
\]
whence, by the counital property, $c = \sum c_i
\epsilon_{\coring{C}}(f(m_i)) \in \coring{J}$. This proves 1.
Statements 2 and 3 are easy consequences of the definition of the
bicomodule of coefficients.
\end{proof}

\section{Semisimple corings}

We study the simplest kind of corings from the categorical point
of view, namely, those corings having a semisimple category of
comodules. We prove generalizations of  known theorems for
coalgebras and rings. In particular, we get a (unique)
decomposition of any semisimple coring in terms of simple
components. The structure of this simple components, which in the
cases of rings and coalgebras over fields is described in terms
of matrices, seems to be much more tangled in the present general
setting. See, however, the last section, for a structure theorem
for the case of simple semisimple corings having a grouplike
element.

\begin{theorem}\label{semi1}
Let $\coring{C}$ be an $A$--coring. The following
statements are equivalent:
\begin{enumerate}[(i)]
\item $\coring{C}$ is semisimple as a left $\coring{C}$--comodule and $\coring{C}_A$ is flat;
\item every left $\coring{C}$--comodule is semisimple and $\coring{C}_A$ is flat;
\item $\coring{C}$ is semisimple as a right $\coring{C}$--comodule and ${}_A\coring{C}$ is flat;
\item every right $\coring{C}$--comodule is semisimple and ${}_A\coring{C}$ is flat;
\item every (left or right) $\coring{C}$--comodule is semisimple, and ${}_A\coring{C}$
and $\coring{C}_A$ are projectives;
\item $\coring{C}$ is semisimple as a left $\coring{C}^*$--module, and as right $\rring$--module,
and ${}_A\coring{C}$ and $\coring{C}_A$ are projectives.
\end{enumerate}
\end{theorem}
\begin{proof}
Since in $(i)$ and $(ii)$ $\coring{C}_A$ is assumed to be flat, we
know that $\rcomod{\coring{C}}$ is a Grothendieck category and,
therefore, the equivalence between $(i)$ and $(ii)$ is a
consequence of Corollary \ref{Csubgenera}. Now, let us show that
$(ii)$ does imply $(vi)$. By Proposition \ref{plano}, the
forgetful functor $U : \lcomod{\coring{C}} \rightarrow \lmod{A}$
is exact. Moreover, it has an exact right adjoint $\coring{C}
\tensor{A} -$, which implies that it preserves projective objects.
Therefore, every left $\coring{C}$--comodule is projective as a
left $A$--module, and, in particular, ${}_A\coring{C}$ is
projective. By Corollary \ref{ratcom}, the category
$\rcomod{\coring{C}}$ of all right $\coring{C}$--comodules is
isomorphic to the category $\Rat(\lmod{{}^*\coring{C}})$ of
rational left ${}^*\coring{C}$--modules. Moreover, since
$\coring{C}$ is a semisimple object in the Grothendieck category
$\lcomod{\coring{C}}$, it follows that $\coring{C}$ is semisimple
as a left module over its endomorphism ring. Now, we know that
this last ring is isomorphic to ${}^*\coring{C}$. Therefore,
${}_{{}^*\coring{C}}\coring{C}$ is semisimple and, thus,
$\coring{C}_{\coring{C}}$ is semisimple.  Our coring satisfies now
conditions symmetric to that in $(ii)$ which entails that
$\coring{C}_{A}$ is also projective and that $\lcomod{\coring{C}}$
is isomorphic to $\Rat(\rmod{\coring{C}^*})$. Of course, we deduce
that $\coring{C}_{\coring{C}^*}$ is semisimple, too, and we arrive
at $(vi)$. Using Corollary \ref{ratcom} we get $(vi) \Rightarrow
(i)$. On the other hand, symmetric arguments are used to show that
$(iii)$, $(iv)$, and $(vi)$ are equivalent. Finally, the
equivalence $(vi) \Leftrightarrow (v)$ is consequence of Corollary
\ref{ratcom}.

%The equivalence between $(vi)$ and $(v)$ and the implications $(v)
%\Rightarrow (iv)$ are consequences of Corollary \ref{ratcom}.
%Finally, $(iv)$ dose obviously imply $(iii)$, and
%$(iii)\Rightarrow (vi)$ follows by symmetry.
\end{proof}
%\item every left $\coring{C}$--comodule is injective $\coring{C}_A$ is flat;
%\item every left $\coring{C}$--comodule is projective $\coring{C}_A$ is flat;
%\item every right $\coring{C}$--comodule is injective ${}_A\coring{C}$ is flat;
%\item every left $\coring{C}$--comodule is projective ${}_A\coring{C}$ is flat;

\begin{definition}
An $A$--coring satisfying the equivalent conditions in Theorem
\ref{semi1} will be called a \emph{semisimple} coring.
\end{definition}

The always marvelous Wedderburn-Artin's structure theorem for
semisimple artinian rings reposes upon a unique decomposition of
the ring as a direct sum of simple artinian rings. This
`abstract' part of that classical result holds in the present
setting. We first define the natural notions of simple coring and
semiartinian coring.

%\begin{corollary}\label{leftright}
%Let $\coring{C}$ be an $A$--coring and assume that the modules ${}_A\coring{C}$ and $\coring{C}_A$
%are projective. Then $\coring{C}$ is left semisimple if and only if  $\coring{C}$ is right
%semisimple.
%\end{corollary}
%\begin{proof}
%Assume that $\coring{C}$ is left semisimple. By Proposition \ref{semi1},
%$\coring{C}_{\coring{C}^*}$ is semisimple.
%This implies that $\coring{C}$ is a semisimple left module over $\rend{\coring{C}^*}{\coring{C}}$.
%By Corollary \ref{ratcom} $\rend{\coring{C}^*}{\coring{C}} = \lend{\coring{C}}{\coring{C}}$ and,
%on the other hand, $\lend{\coring{C}}{\coring{C}} \cong {}^*\coring{C}$ as rings. Therefore,
%${}_{\coring{C}^*}\coring{C}$ is semisimple. The symmetric version of Proposition \ref{semi1}
%gives that $\coring{C}$ is right semisimple.
%\end{proof}

\begin{definition}
A coring is said to be \emph{simple} if it does not contain
non-trivial sub-bicomodules. Notice that if $\coring{C}$ is a
semisimple coring, then it is simple if and only if it does not
contain non-trivial sub-corings.
\end{definition}

\begin{definition}
Assume that the category of all left comodules over an
$A$--coring $\coring{C}$ is a Grothendieck category (see
Proposition \ref{plano}). The coring $\coring{C}$ is said to be
left \emph{semiartinian} if it is semiartinian as an object in
$\lcomod{\coring{C}}$, namely, every factor comodule of
${}_{\coring{C}}\coring{C}$ contains a (nonzero) simple
subcomodule.
\end{definition}

\begin{theorem}\label{simpsemisimp}
Let $\coring{C}$ be an $A$--coring such that the modules ${}_A\coring{C}$ and $\coring{C}_A$
are projective. The following statements are equivalent:
\begin{enumerate}[(i)]
\item $\coring{C}$ is a simple left semi-artinian coring;
\item $\coring{C}$ is a simple coring and contains a simple left $\coring{C}$--subcomodule;
\item $\coring{C}$ is a semisimple coring with a unique type of simple left
$\coring{C}$-comodule;
\item $\coring{C}$ is a simple right semi-artinian coring;
\item $\coring{C}$ is a simple coring and contains a simple right $\coring{C}$--subcomodule;
\item $\coring{C}$ is a semisimple coring with a unique type of simple right
$\coring{C}$-comodule.
\end{enumerate}
\end{theorem}
\begin{proof}
$(i) \Rightarrow (ii)$ is obvious.\\
$(ii) \Rightarrow (iii)$ Let $S$ be a simple left subcomodule of $\coring{C}$. By Corollary
\ref{ratcom}, $S$ is a simple right $\coring{C}^*$--submodule of $\coring{C}$. Now,
${}^*\coring{C}S$ is a nonzero $({}^*\coring{C},\coring{C}^*)$--bi-submodule of $\coring{C}$,
which, by Corollary \ref{subcosubbi}, is a nonzero sub-bicomodule.
Hence, ${}^*\coring{C}S = \coring{C}$ and, therefore, $\coring{C}$ is a sum of homomorphic
images of the simple right $\coring{C}^*$--module $S$. Apply Corollary \ref{ratcom}.\\
$(iii) \Rightarrow (i)$ Obviously, every semisimple coring is
left semi-artinian. Let $\coring{I}$ be a non-zero sub-bicomodule
of $\coring{C}$. In particular, $\coring{I}$ is a left
$\coring{C}$--subcomodule of $\coring{C}$, so that it contains a
simple subcomodule $S$. By the statements 1 and 2 of Proposition
\ref{C(M)} we get that $C(S) \subseteq C(\coring{I}) =
\coring{I}$. Since $\coring{C}$ is isomorphic to a direct sum of
copies of $S$, we apply part 3 of Proposition \ref{C(M)} to
obtain $C(S) = \coring{C}$.
Hence, $\coring{I} = \coring{C}$ and $\coring{C}$ is simple. \\
$(ii) \Rightarrow (v)$ We have already proved that if
$\coring{C}$ is simple and contains a simple left comodule, then
$\coring{C}$ is semisimple. Thus, it contains a simple
right $\coring{C}$--comodule. \\
Finally, $(iv), (v), (vi)$ are proved to be equivalent in a analogous way to the
proof of the equivalence between $(i), (ii), (iii)$;
which allows to derive also that $(v)$ implies $(ii)$. This finishes the proof.
\end{proof}

\begin{remark}\label{intrincado} Let $\coring{C}$ be a simple semi-artinian
$A$--coring. By Theorem \ref{simpsemisimp} we have that
${}_{\coring{C}}\coring{C} \cong S^{(\Xi)}$, where $S$ is a simple
left $\coring{C}$--comodule and $\Xi$ is an index set. In contrast
with the coalgebra or ring cases (i.e., when $\coring{C}$ is a
coalgebra over a field or $\coring{C} = A$), the set $\Xi$ needs
not to be finite. In fact, consider the $A$--coring $I$ given in
Example \ref{elejemplo} with $R$ a simple artinian ring, $B$ the
coproduct of $\Xi$ copies of the unique simple left $R$--module
($\Xi$ is any infinite set) and $S$ is the endomorphism ring of
the left $R$--module $B$. Then $I$ is a simple semi-artinian ring
and it is isomorphic to a direct sum of infinitely many copies of
a simple left $I$--comodule (which is essentially the unique
simple left $R$--module). This easy example also shows that the
basis ring $A$ needs not to be semisimple or even artinian for a
semisimple $A$--coring.
\end{remark}

We finish this section by showing that semisimple corings can be completely described
in terms of simple semiartinian (or simple semisimple corings).

\begin{theorem}\label{semisimple} An $A$--coring $\coring{C}$ is semisimple if and only if
it decomposes as $\coring{C} = \oplus_{\alpha \in
\Lambda}\coring{C}_{\alpha}$, where $\coring{C}_{\alpha}$ is a
simple semisimple $A$--subcoring for every $\alpha \in \Lambda$.
In such a case, the decomposition is unique.
\end{theorem}
\begin{proof} Assume that $\coring{C}$ is semisimple. Let $\Lambda$ be a set of
representatives of all simple right $\coring{C}$--comodules. For
each $\alpha \in \Lambda$, define $\coring{C}_{\alpha}$ to be the
$\alpha$--th isotypic component of $\coring{C}_{\coring{C}}$.
Since $\coring{C}$ is right semisimple, it follows that
$\coring{C} = \oplus_{\alpha \in \Lambda}\coring{C}_{\alpha}$. We
know from Corollary \ref{ratcom} that $\coring{C}_{\alpha}$ is a
left ${}^*\coring{C}$--submodule of $\coring{C}$. Given $c^* \in
\coring{C}^*$, its right multiplication map is a homomorphism of
left ${}^*\coring{C}$--modules, and, thus, of right
$\coring{C}$--comodules. It follows that $\coring{C}_{\alpha}$ is
a right $\coring{C}^*$--submodule of $\coring{C}$ and, by
Corollary \ref{subcosubbi}, $\coring{C}_{\alpha}$ is a subcoring
of $\coring{C}$. Obviously, $\coring{C}_{\alpha}$ is semisimple
with a unique type of simple; by Theorem \ref{simpsemisimp},
$\coring{C}_{\alpha}$ is a simple semi-artinian $A$--coring.
Finally, the converse implication is easily deduced from the fact
that, given the stated decomposition $\coring{C} = \oplus_{\alpha
\in \Lambda}\coring{C}_{\alpha}$, the right
$\coring{C}$--subcomodules of $\coring{C}$ are of the form
$\oplus_{\alpha \in \Lambda} M_{\alpha}$, where $M_{\alpha}$ is a
$\coring{C}_{\alpha}$--subcomodule of $\coring{C}_{\alpha}$ for
every $\alpha$. The uniqueness comes from the observation that the
$\coring{C}_{\alpha}$'s are just the isotypic components of
$\coring{C}_{\coring{C}}$.
\end{proof}

\section{ Simple semiartinian corings with a grouplike element}

A complete description of all semisimple corings over a given ring $A$ would be obtained,
in view of Theorem \ref{semisimple}, throughout the knowledge of the structure of simple
semiartinian $A$--corings. The structure of a general simple semiartinian coring
seems to be quite intricate (see Example \ref{intrincado}). It is possible, however, to recognize
the simple semiartinian $A$--corings having a grouplike element as the canonical corings
$A \tensor{B} A$, where $B$ runs the set of simple artinian subrings of $A$, as we will
prove in this section.

Let $\coring{C}$ be an $A$--coring, a non-zero element $g \in
\coring{C}$ such that $\epsilon(g) =1$ and $\Delta(g) = g
\tensor{A} g$ is called a \emph{grouplike} element. An example of
coring with such an element is $A \tensor{B}A$ cited in Example
\ref{ejemplo1} taking $g=1 \tensor{B} 1$.

\begin{lemma}\label{Acomod}\cite[Lemma 5.1]{Brzezinski:2000unp}
Let $\coring{C}$ be an $A$--coring. Then $A$ is a right
$\coring{C}$--comodule if and only if $A$ is a left
$\coring{C}$--comodule if and only if there exists a grouplike
element $g \in \coring{C}$. In that case the left and right
coactions are given by
\[
\xymatrix@R=0pt{ \lambda_A : A \ar@{->}[r] & \coring{C} \\ a
\ar@{|->}[r] & a g\tensor{A} 1 }, \xymatrix@R=0pt{\rho_A : A
\ar@{->}[r] & \coring{C} \\ a \ar@{|->}[r] & 1\tensor{A}ga. }
\]
\end{lemma}

Assume that $\coring{C}$ has a grouplike element $g$, and consider
the \emph{subring of coinvariants} of $A$ defined by
\[
A^{co\coring{C}}=\{a \in A|\, ag=ga \};
\]
this ring is isomorphic to $\rend{\coring{C}}{A}$, and also to
$\lend{\coring{C}}{A}$. Then we have a functor \cite[Proposition
5.2]{Brzezinski:2000unp} $(-)^{co\coring{C}} : \rcomod{\coring{C}}
\rightarrow \rmod{A^{co\coring{C}}}$ which assigns to every right
$\coring{C}$--comodule $M$ the right $A^{co\coring{C}}$--module of
\emph{coinvariants}
\[
M^{co\coring{C}} = \{ m \in M | \rho_M (m) = m \tensor{A} g \}.
\]
It is easily shown that this functor is naturally isomorphic to
the functor $\hom{\coring{C}}{A_{\coring{C}}}{-}$. The analogous
discussion is pertinent for left $\coring{C}$--comodules.

\begin{proposition}\label{simple-artin} Let $ B \rightarrow A$ be a
ring extension, and $A \tensor{B} A$ with the canonical
$A$--coring structure defined in Example \ref{ejemplo1}. Assume
that $B$ is a simple artinian ring. Then $A \tensor{B} A$ is a
simple semisimple $A$--coring. Moreover, $A^{co A \tensor{B} A} =
B$, with respect to the grouplike element $1 \tensor{B}1$.
\end{proposition}
\begin{proof} We
know that $A_B$ and ${}_BA$ are projective modules, and this
implies that the coring $\coring{C} = A \tensor{B} A$ is
projective as a left and as a right $A$--module. By Corollary
\ref{ratcom}, the category $\lcomod{A \tensor{B} A}$ is isomorphic
to the category of all rational right $\coring{C}^*$--modules.
Recall from \cite[Example 3.3]{Sweedler:1975} that there is an
anti-isomorphism of rings \[ \xymatrix@R=0pt{ \rring \ar@{->}[r] &
\End{A_B} \\ g \ar@{|->}[r] & [ a \mapsto g(a \tensor{B} 1)] \\ f
\tensor{B} A  \ar@{<-|}[r] & f. } \] Some straightforward
computations show that the canonical left $\End{A_B}$--module
structure of $A$ corresponds to the structure of (rational) right
$\coring{C}^*$--module, whenever the coactions of $A$ are derived
from $1 \tensor{B}1$ (see Lemma \ref{Acomod}). Since $A_B$ is a
homogeneous semisimple right $B$--module it follows that ${
}_{\End{A_B}}A$ is homogeneous semisimple, too. We conclude, by
Lemma \ref{Acomod} and Theorem \ref{thm:iso}, that $A$ is direct
sum of copies of a simple left $\coring{C}$--comodule. Now, let $a
\in A$, and consider the following homomorphism of left
$A$--modules \[ \xymatrix@R=0pt{ \phi_a : A \ar@{->}[r] &
\coring{C} \\ a' \ar@{|->}[r] & a' \tensor{B} a, } \] which is, in
fact, a homomorphism of left $\coring{C}$--comodules. It follows
that $A$ generates $A \tensor{B} A$ as a left comodule. In
particular, we get that $A \tensor{B} A$ is a sum of copies of a
simple comodule which, in the light of Theorem \ref{simpsemisimp},
shows that $A \tensor{B} A$ is a simple semisimple $A$--coring.
Finally,
\[
A^{co A \tensor{B} A} \cong \End{{}_{A \tensor{B} A}A} =
\End{A_{\coring{C}^*}} = \End{{}_{\End{A_B}}A} =
\mathrm{Biend}(A_B),
\]
and $A_B$ is a balanced $B$--module (remember that $B$ is simple
artinian). Hence, $B = A^{co A \tensor{B} A}$.
\end{proof}

Next theorem tell us that Proposition \ref{simple-artin} gives
all possible examples of simple semisimple corings with a
grouplike element.

\begin{theorem}\label{invariantes}
Let $\coring{C}$ be an $A$--coring, and $g \in \coring{C}$ be a
grouplike element. Assume that $\coring{C}$ is a simple semisimple
$A$--coring. Then $A^{co \coring{C}}$ is a simple artinian ring
and the canonical $A$--bimodule map $A \tensor{A^{co \coring{C}}}
A \rightarrow \coring{C}$ which sends $1 \tensor{A^{co
\coring{C}}} 1$ to $g$ is an isomorphism of $A$--corings.
\end{theorem}
\begin{proof}
Endow $A$ with the structure of right $\coring{C}$--comodule
derived from $g$. Since $\coring{C}$ is assumed to be simple
semisimple, $A$ is a direct sum of copies of the only simple
right $\coring{C}$--comodule. Moreover, this direct sum, being of
right $A$--submodules after all, is finite. Therefore,
$\End{A_{\coring{C}}} \cong A^{co\coring{C}}$ is a simple artinian
ring.
 It is easily proved that the $A$--bimodule homomorphism
 $\varphi: A \tensor{A^{\coring{C}}} A \rightarrow \coring{C}$ which sends $a \tensor{} a'$ onto
 $aga'$ is an $A$--coring homomorphism. Since $A_{\coring{C}}$ is
 a finitely generated projective generator for
 $\rcomod{\coring{C}}$, a standard consequence of Gabriel-Popescu's Theorem (see, e.g.
 \cite[Corolar 9.7]{Nastasescu:1976}) says that $\hom{\coring{C}}{A}{-} :
 \rcomod{\coring{C}} \rightarrow \rmod{\rend{\coring{C}}{A}}$ is an equivalence of categories.
 Now, $\hom{\coring{C}}{A}{-} \cong
 (-)^{co\coring{C}}$ naturally which implies, by \cite[Theorem
 5.6]{Brzezinski:2000unp}, that $\varphi$ is an isomorphism.
\end{proof}

Following \cite[Definition 5.3]{Brzezinski:2000unp} we say that an
$A$--coring with grouplike $g$ is \emph{Galois} if the $A$--coring
map which sends $1 \tensor{A^{co\coring{C}}} 1$ to $g$ gives an
isomorphism $\coring{C} \cong A \tensor{A^{co\coring{C}}} A$.
Thus, Theorem \ref{invariantes} says that every simple semisimple
$A$--coring with a grouplike element is Galois. We have already
more, as the following theorem, which collects the relevant
information about simple semisimple corings with a grouplike
element, shows.

\begin{theorem}
The following conditions are equivalent for an $A$--coring
$\coring{C}$ with a grouplike element $g$.
\begin{enumerate}[(i)]
\item $\coring{C}$ is a simple simisimple $A$--coring;
\item $\coring{C} \cong A \tensor{B} A$ for some simple artinian
subring $B$ of $A$;
\item $\coring{C}$ is Galois and $A^{co\coring{C}}$ is a simple
artinian ring;
\item $\coring{C}_A$ is flat, $A$ is a projective generator in
$\lcomod{\coring{C}}$, and $A^{co\coring{C}}$ is a simple artinian
ring;
\item ${}_{A}\coring{C}$ is flat, $A$ is a projective generator in
$\rcomod{\coring{C}}$, and $A^{co \coring{C}}$ is a simple
artinian ring.
\end{enumerate}
\begin{proof}
$(i) \Rightarrow (iii)$ This is Theorem \ref{invariantes}. \\
$(iii) \Rightarrow (ii)$ Obvious. \\
$(ii) \Rightarrow (i)$ It follows from Theorem \ref{simple-artin}. \\
$(i) \Rightarrow (iv)$ By Theorem \ref{semi1}, $\coring{C}_A$ is in fact projective. Theorem
\ref{simpsemisimp} gives that every right left $\coring{C}$--comodule is a direct sum of copies
of the unique simple comodule. Thus, every nonzero comodule is a projective generator for
$\lcomod{\coring{C}}$. Finally, since $(i)$ is equivalent to $(iii)$, we know that
$A^{co\coring{C}}$ is a simple artinian ring.\\
$(iv) \Rightarrow (iii)$ The proof of Theorem \ref{invariantes} is easily adapted to obtain this
implication.\\
$(v) \Leftrightarrow (iv)$ It follows by symmetry.
\end{proof}
\end{theorem}

%\bibliographystyle{amsplain}
%\bibliography{c:/texdoc/bib/nueva}

\providecommand{\bysame}{\leavevmode\hbox
to3em{\hrulefill}\thinspace}
\providecommand{\MR}{\relax\ifhmode\unskip\space\fi MR }
\providecommand{\MRhref}[2]{%
  \href{http://www.ams.org/mathscinet-getitem?mr=#1}{#2}
} \providecommand{\href}[2]{#2}

\begin{quote}
 \small{\textsc{Departamento de Algebra, Facultad de
Ciencias, Universidad de Granada, E18071-Granada, Spain.}}
\end{quote}

\end{document}